\documentclass[12pt]{article}
\usepackage{graphicx}
\usepackage{tikz}
\usepackage{cite}
\usepackage{hyperref}
\usepackage{amsfonts}
\usepackage{fancyhdr,amsthm,amsmath,amssymb,mathrsfs}
\usepackage{bookmark}
\usetikzlibrary{intersections}

\newtheorem{thm}{\bf Theorem}[section]
\newtheorem{lem}[thm]{Lemma}
\newtheorem{prop}[thm]{Proposition}
\newtheorem{cor}[thm]{Corollary}
\theoremstyle{definition}

\numberwithin{equation}{section}
\newenvironment{prf}{\noindent {\bf Proof.}\rm}{\qed}

\newcommand{\DOI}[1]{DOI: \href{https://doi.org/#1}{#1}}

\begin{document}
\title{\Large\bf Distributivity in lattices of torsion classes over finite-dimensional algebras\thanks{Supported by the National Natural Science
Foundation of China (Grant No. 11771212) and the Postgraduate Research $\&$ Practice Innovation Program of Jiangsu Province (No. KYCX22-1533).}}
\author{Yongle Luo and Jiaqun Wei\thanks{Corresponding author} \\
{\small\em School of Mathematical Sciences, Nanjing Normal University, $210023$ Nanjing, China} \\
{\small\em E-mail: yongleluo@nnu.edu.cn, \;  weijiaqun@njnu.edu.cn} }
\date{}
\maketitle

\begin{abstract}
Let $A$ be a basic finite-dimensional algebra and denote by $\operatorname{tors} A$ the collection of all all torsion classes of $A$. It has been proved in \cite{Demonet} that $\operatorname{tors} A$ is always a completely semidistributive lattice. In the present paper, we investigate the distributivity of this lattice and proved that the lattice is distributive if and only if the algebra is the finite direct product of finite-dimensional local algebras.
\medskip


\end{abstract}

\section{Introduction and main result}
Torsion theory is an important concept in representation theory. It is closely related to the study of derived categories and their $t$-structures. Recently, Demonet-Iyama-Reading-Reiten-Thomas in \cite{Demonet} established a lattice theoretical framework to study the torsion classes over a finite-dimensional algebra. Based on these descriptions, we considered the distributivity of both the lattice of torsion classes and the poset of functorially finite torsion classes. Then we proved that a basic finite-dimensional algebra with this property is precisely the finite product of local algebras.

Let $A$ be a finite-dimensional algebra over an algebraically closed field $K$. We denote by $\operatorname{mod} A$ the category of finite-dimensional right $A$-modules. Throughout this paper, all subcategories of $\operatorname{mod} A$ are full. A subcategory of $\operatorname{mod} A$ is called a {\it torsion class} (respectively, {\it torsion-free class}) if it is closed under extensions and factor modules (respectively, submodules). Given a subcategory $\mathscr{X}$ of  $\operatorname{mod} A$, we define its {\it $\operatorname{Hom}$-orthogonal subcategories} of $\operatorname{mod} A$ by
\[\mathscr{X}^\perp= \{M\in \operatorname{mod} A| \mbox{Hom}_{A}(X, M)=0, \forall X\in \mathscr{X}\},
\]
\[^\perp\mathscr{X}= \{N\in \operatorname{mod} A| \mbox{Hom}_{A}(N, X)=0, \forall X\in \mathscr{X}\}.
\]

Let $\mathscr{T}$ and $\mathscr{F}$ be two subcategories of $\operatorname{mod} A$. We call $(\mathscr{T}, \mathscr{F})$ a {\it torsion pair} if $\mathscr{T}^\perp = \mathscr{F}$ and $^{\perp}\mathscr{F} = \mathscr{T}$.
For a torsion pair $(\mathscr{T}, \mathscr{F})$, it is well known that $\mathscr{T}$ is a torsion class uniquely determined by $\mathscr{F}$, and $\mathscr{F}$ is a torsion-free class uniquely determined by $\mathscr{T}$.

Denote by $\operatorname{tors} A$ (respectively, $\operatorname{torf} A$) the set of all torsion classes (respectively, torsion-free classes). Iyama-Reiten-Thomas-Todorov \cite{Iyama} proved that the two posets $\operatorname{tors} A$ and $\operatorname{torf} A$ ordered by inclusion are both complete lattices. In addition, the two complete lattices are anti-isomorphic:
\[(-)^{\perp}: \operatorname{tors} A \rightarrow \operatorname{torf} A, ~~~^{\perp}(-): \operatorname{torf} A \rightarrow \operatorname{tors} A.
\]

A module $B \in \operatorname{mod} A$ is called a {\it brick} if every non-zero endomorphism of $B$ is an isomorphism. Barnard-Carroll-Zhu \cite{Barnard} proved that the {\it completely join irreducible elements} of $\operatorname{tors} A$ are in bijection with the {\it bricks} of $\operatorname{mod }A$. By considering the Hasse quiver of the lattice of torsion classes, Demonet-Iyama-Reading-Reiten-Thomas \cite{Demonet} used the {\it brick labeling} (see also in \cite{Asai,Barnard}) and showed that $\operatorname{tors} A$ is completely semidistributive (see also in \cite{Garver}), bialgebraic and hence arrow-separated. For more details of the semidistributivity of this lattice, one can also refer to \cite{Thomas}.

Since the poset of torsion classes forms a complete semidistributive lattice, an essential question is whether this lattice satisfies the distributive property, which plays a crucial role in combinatorial and structural representation theory. A lattice $L$ is called {\it distributive} if the meet operation distributes over the join and vice-versa. For instance, given a poset $(P,\leq)$ we consider $\operatorname{Ideal}(P)$, the set of {\it order ideals} of $P$ (downwards closed subsets). The inclusion of ideals makes $\operatorname{Ideal}(P)$ into a distributive lattice, where the meet is given by intersection of ideals and the join is given by union of ideals. Recall from {\it Birkhoff representation theorem}, any finite distributive lattice $L$ is isomorphic to $\operatorname{Ideal}(P)$ where $P$ is the subposet of $L$ containing the {\it join-irreducible elements} of $L$.

The main purpose of the present paper is to study the distributivity of $\operatorname{tors} A$. It is easy to see that the lattice of torsion classes over a semisimple algebra is distributive. However, determining which algebras have this property remains a nontrivial problem. In this work, we show that the lattice is distributive if and only if every torsion class is completely determined by its simple modules, meaning that the set of bricks coincides with the set of simple modules. Moreover, in the case of a basic finite-dimensional algebra, this condition implies that the algebra is precisely the finite direct product of finite-dimensional local algebras. Our main result can be stated as follows.
\begin{thm}\label{Main theorem}
Let $A$ be a basic finite-dimensional algebra over an algebraically closed field $K$.  Then the lattice of torsion classes $\operatorname{tors}A$ is distributive
if and only if $A$ is the finite direct product of finite-dimensional local algebras.
\end{thm}

Recall from \cite{smalo} that a torsion class $\mathscr{T}$ in $\operatorname{mod} A$ is {\it functorially  finite} precisely when it is of the form $\operatorname{Fac}(M)$ for some $M \in \operatorname{mod} A$, where $\operatorname{Fac}(M)$ is the full subcategory of $\operatorname{mod} A$ consisting of all factor modules of finite direct sums of copies of $M$. The functorially finite torsion class plays a fundamental role in the study of torsion theory, as it provides a natural bridge to $\tau$-tilting theory. In particular, Adachi-Iyama-Reiten\cite{Adachi} established a bijection between functorially finite torsion classes and the isomorphism classes of basic support $\tau$-tilting modules. This connection is crucial in modern representation theory, as $\tau$-tilting theory generalizes classical tilting theory and provides deep insights into the structure of module categories.

Denote by $\operatorname{f-tors} A$ the poset of all functorially finite torsion classes in $\operatorname{tors} A$ ordered by inclusion. In general, $\operatorname{f-tors} A$ is a subposet of $\operatorname{tors} A$, but it need not form a lattice. Iyama-Reiten-Thomas-Todorov \cite{Iyama} considered the case when $A=KQ$, where $KQ$ is the path algebra of a finite connected acyclic quiver $Q$. They showed that $\operatorname{f-tors} KQ$ forms a lattice (not necessarily complete), if and only if $Q$ is either a Dynkin quiver or has at most two vertices. In these cases, Yang \cite{Yang} proved that the set of isomorphism classes of basic support $\tau$-tilting $KQ$-modules ($\operatorname{f-tors} KQ$) forms a distributive lattice if and only if $Q$ is of type $A_1$.

In the present paper, we also investigate the distributivity of $\operatorname{f-tors} A$ for any finite dimensional $K$-algebra $A$. And we show that $\operatorname{f-tors} A$ is a distributive lattice only if $A$ is $\tau$-tilted finite, leading to $\operatorname{f-tors} A= \operatorname{tors} A$. Then by Theorem \ref{Main theorem}, we have the following corollary which extends the result of \cite[Theorem 1.4]{Yang}.

\begin{cor} Let $A$ be a basic finite-dimensional $K$-algebra. Then $\operatorname{f-tors} A$ is a distributive lattice if and only if $A$ is the finite product of finite-dimensional local algebras.
\end{cor}

\section{Preliminaries}
\subsection{Lattice theory}
We recall some fundamental materials on lattices that can be found in the standard lattice theory book \cite{Gratzer}.

Given a poset $L$ and $a, b \in L$, we say that $a$ {\it covers} $b$, and write $a\succ b$, if $a > b$ and there is no element $c$ such that $a > c > b$. The {\it Hasse quiver} of a poset $L$ is the directed graph whose vertices correspond to the elements of $L$, and there is an arrow from $a$ to $b$ if and only if $a \succ b$.  Two elements $a$ and $b$ are {\it adjacent} in a Hasse quiver, if there is an arrow from $a$ to $b$ ($a\succ b$), or an arrow from $b$ to $a$ ($b \succ a$). A poset $L$ is called {\it Hasse-regular} if each vertex of its Hasse quiver has the same degree when considered as an undirected graph.

A poset $L$ is called a {\it lattice} if, for all $a, b \in L$, there exists a least upper bound $a\vee b\in L$, called the {\it join} of $a$ and $b$, and a greatest lower bound $a\wedge b \in L$, called the {\it meet} of $a$ and $b$. Let $L=(L,\leq, \vee, \wedge)$ be a lattice. Recall that $L$ is {\it upper semimodular} if $a, b\succ a\wedge b$ implies that $a \vee b \succ a, b$, for all $a, b \in L$. Dually, $L$ is {\it lower semimodular} if
$a\vee b \succ a, b$ implies that $a, b\succ a\wedge b$, for all $a, b \in L$. Moreover, $L$ is {\it distributive} if $(a\vee b)\wedge c =(a\wedge c)\vee (b\wedge c)$, for all $a, b, c \in L$. In particular, a bounded lattice $L$ (i.e., $L$ has a minimum 0 and a maximum 1) is
{\it Boolean} if $L$ is distributive and for any $a\in L$, there exists $b\in L$ such that $a\vee b= 1$ and $a\wedge b=0$. For instance, the class of all subsets of $\{1, \dots, n\}$  ordered by inclusion forms a Boolean lattice.

A lattice $L$ is {\it complete} if any subset $S$ of $L$ has a join, which we denote $\bigvee\nolimits_{x\in S} x$ or $\bigvee S$, and a meet, which we denote $\bigwedge\nolimits_{x\in S} x$ or $\bigwedge S$. Given a complete lattice $L$ and $x \in L$, $x$ is called {\it completely join irreducible} if $\bigvee\nolimits_{y < x} y < x$.

The complete lattice $L$ is {\it completely join semidistributive}, if given $x \in L$ and a set $S \subseteq L$ such that $x \vee y = z$ for all $y \in S$, then $x \vee (\bigwedge S) =z$; The definition of {\it completely meet semidistributive} is defined dually. The complete lattice $L$ is {\it completely  semidistributive}, if $L$ is both completely join semidistributive and completely meet semidistributive.

\subsection{Lattices of $\mathbf{tors} A$ and $\mathbf{torf} A$}
Each torsion class (respectively, torsion-free class) can be determined by filtration of modules. Let $\mathcal{C}$ be a subcategory of $\operatorname{mod}A$. We define $\mathscr{T}(\mathcal{C})$ (respectively, $\mathscr{F}(\mathcal{C})$) the subcategory whose modules are filtered by quotients (respectively, submodules) of modules from $\mathcal{C}$. It is a classical result that $\mathscr{T}(\mathcal{C})$ (respectively, $\mathscr{F}(\mathcal{C})$) is the smallest torsion class (respectively, torsion-free class) containing all modules from $\mathcal{C}$ (see \cite[Proposition 2.1]{Thomas}). Then it follows that $\mathscr{T}(\mathcal{C})$ = $^{\bot}(\mathcal{C}^{\perp})$ with induced torsion pair ($\mathscr{T}(\mathcal{C}), \mathcal{C}^{\perp})$, and $\mathscr{F}(\mathcal{C})= (^{\perp}\mathcal{C})^{\perp}$ with induced torsion pair ($^{\perp}\mathcal{C}, \mathscr{F}(\mathcal{C}))$.

Since the intersection of torsion classes is again a torsion class. The poset $\operatorname{tors}A$ ordered by inclusion clearly has a meet given by intersection. Dually, the same is true for $\operatorname{torf} A$. Note that the left (respectively, right) orthogonal operation is order-reversing bijections between $\operatorname{tors} A$ and $\operatorname{torf} A$. Let $S$ be any subset of $\operatorname{tors} A$, then we have
\[\bigvee_{\mathcal{T} \in S} \mathcal{T}^\perp
= {}^{\perp\!} \bigl( \bigwedge_{\mathcal{T} \in S} \mathcal{T}^\perp \bigr).
\]

Since the join is the smallest torsion class which contains all $\mathcal{T}$ in $S$, we can also define the join directly as follows:
\[\bigvee_{\mathcal{T}\in S} \mathcal{T} = \mathscr{T}(\bigcup_{\mathcal{T}\in S} \mathcal{T}).
\]

\subsection{$\tau$-tilting theory}
Recall that a module $T\in \operatorname{mod} A$ is {\it $\tau$-rigid} if $\operatorname{Hom}_{A}(T, \tau T)=0$, where $\tau$ is the Auslander-Reiten translation. A module $T\in \operatorname{mod} A$ is called {\it $\tau$-tilting} if it is $\tau$-rigid and has $n$ non-isomorphic indecomposable summands, where $n$ is the number of non-isomorphic simple $A$-modules. Moreover, we say that $T\in \operatorname{mod} A$ is {\it support $\tau$-tilting} if it is a $\tau$-tilting $A/\langle e \rangle$-module for some idempotent $e\in A$. Denote by $\operatorname{s\tau-tilt} A$
(respectively, $\operatorname{i\tau-rigid} A$) the set of isomorphism classes of basic support $\tau$-tilting $A$-modules (respectively, indecomposable $\tau$-rigid $A$-modules).

Let $M, N \in \operatorname{s\tau-tilt} A$, we say that $M$ and $N$ are {\it mutations} of each other if there is an  $A$-module $U$ such that $M=U\oplus X$, $N=U\oplus Y$ and $X$, $Y$ are either $0$ or indecomposable. By \cite[Definition-Proposition 2.28]{Adachi}, either $\operatorname{Fac}(N)\subsetneq \operatorname{Fac}(M)$ or $\operatorname{Fac}(M)\subsetneq \operatorname{Fac}(N)$ holds. In this case, then we say that $N$ is a {\it left mutation} (respectively, {\it right mutation}) of $M$ if $\operatorname{Fac}(N) \subsetneq \operatorname{Fac}(M)$ (respectively, $\operatorname{Fac}(M) \subsetneq \operatorname{Fac}(N)$).

The Hasse quiver of the poset $\operatorname{f-tors} A$ ordered by inclusion is completely determined by the mutations of $ \operatorname{s\tau-tilt} A$. Indeed, let $Q(\operatorname{s\tau-tilt} A)$ be the quiver that the set of vertices is $ \operatorname{s\tau-tilt} A$, and draw an arrow from $M$ to $N$ if $N$ is a left mutation of $M$. The Hasse quiver of $\operatorname{s\tau-tilt} A$ is isomorphic to $Q(\operatorname{s\tau-tilt} A)$ (see \cite[Corollary 2.34]{Adachi}).

The algebra $A$ is called {\it $\tau$-tilting finite} if there are only finitely many indecomposable $\tau$-rigid $A$-modules. It follows from \cite [Theorem 1.2]{Iyama} that $A$ is $\tau$-tilting finite if and only if $\operatorname{tors}A= \operatorname{f-tors} A$.

\begin{lem}\label{Hasse regular of ftorsA}\cite[Corrollary 4.6]{Demonet}
Let $A$ be a finite-dimensional $K$-algebra. Then the poset $\operatorname{f-tors} A$ is Hasse-regular and the degree of any vertex is $n$, where $n$ is the number of non-isomorphic simple $A$-modules.
\end{lem}

\section{Proof of Main Results}

The following meaningful results help us describe the algebra $A$ for which $\operatorname{tors} A$ is distributive. Denote by $\operatorname{brick} A$ the set of all bricks of $\operatorname{mod} A$.

\begin{lem}\label{bricks=join irre}\cite[Theorem 1.5]{Barnard}
The map $B \rightarrow \mathscr{T}(B)$ is a bijection from $\operatorname{brick} A$ to the set of all completely join irreducible elements in $\operatorname{tors}A$.
\end{lem}

\begin{lem}\label{every torsion and its bricks}\cite[Lemma 7.1]{Thomas}
Let $\mathcal{T}\in \operatorname{tors} A$. Then $\mathcal{T}$ is characterized by the bricks it contains, that is,
\[\mathcal{T}= \bigvee_{B\in \mathcal{T}\cap \operatorname{brick} A} \mathscr{T}(B).
\]
\end{lem}

\begin{prop}\label{equivalent conditions of torsA}
Let $A$ be a finite-dimensional $K$-algebra. Then the following statements are equivalent.
\begin{description}
\item $(a)$~The lattice $\operatorname{tors} A$ is upper semimodular.
\item $(b)$~Every brick module in $\operatorname{mod} A$ is simple.
\item $(c)$~The lattice $\operatorname{tors} A$ is isomorphic to the Boolean lattice of $\{1, \dots , n\}$ ordered by inclusion of subsets, where $n$ is the number of non-isomorphic simple $A$-modules.
\end{description}
\end{prop}

\begin{prf}
It is clear that $(c)$ implies $(a)$.

$(a)\Rightarrow (b)$: Without loss of generality, we assume that the algebra $A$ has $n$ non-isomorphic simple $A$-modules: $S_1,\dots, S_n$. Note that the arrows ending at $0$ in the Hasse quiver of $\operatorname{tors}A$ are $\mathscr{T}(S_i)\rightarrow 0$, for $i= 1, \dots, n$. Since $S_i$ can not be filtered by the other simple $A$-modules, it is clear that
\[\mathscr{T}(\{S_i\}_{i\in I}) = \mathscr{T}(\{S_j\}_{j\in J}) \Longleftrightarrow I = J,
\]
for any subsets $I, J$ of $\{1, \dots, n\}$. Moreover, we also have
\[\mathscr{T}(\{S_i\}_{i\in I})\bigvee \mathscr{T}(\{S_j\}_{j\in J})= \mathscr{T}(\{S_k\}_{k\in I\cup J}),
\]
\[\mathscr{T}(\{S_i\}_{i\in I})\bigwedge \mathscr{T}(\{S_j\}_{j\in J}) = \mathscr{T}(\{S_k\}_{k\in I\cap J}).\]
Particularly, $\mathscr{T}(\emptyset)=0$ and $\mathscr{T}(\{S_k\}_{k\in \{1, \dots, n\}}) = \operatorname{mod} A$.

Given two non-isomorphic simple $A$-modules $S_i$, $S_j$, it is clear that they are covered by their join, as $\operatorname{tors}A$ is upper semimodular. That is, $\mathscr{T}(S_{i}, S_{j}) \succ \mathscr{T}(S_{i})$, and $\mathscr{T}(S_{i}, S_{j}) \succ\mathscr{T}(S_{j})$. By induction, we get that $$\mathscr{T}(\{S_i\}_{i\in I})\succ \mathscr{T}(\{S_j\}_{j\in J})\Longleftrightarrow I\succ J,$$
for any subsets $I, J$ of $\{1, \dots, n\}$. We remind the reader that the notation $I\succ J$, as introduced earlier, means that $I$ covers $J$, that is, there is no subset $K$ of $\{1, \dots, n\}$ such that $J \subsetneq K \subsetneq I$.

Let $\mathcal{T}$ be an arbitrary torsion class in $\operatorname{tors}A$ and $S_{\mathcal{T}}$ be the set of all simple $A$-modules in $\mathcal{T}$. We claim that $\mathcal{T} = \mathscr{T}(S_{\mathcal{T}})$. Firstly, $\mathscr{T}(S_{\mathcal{T}}) \subseteq \mathcal{T}$. Suppose that $\mathscr{T}(S_{\mathcal{T}}) \subsetneq \mathcal{T}$. Then $\mathscr{T}(S_{\mathcal{T}}\cup\{S_i\})$ is not comparable with $\mathcal{T}$, for any simple $A$-module $S_i \notin S_{\mathcal{T}}$. Indeed, if $\mathscr{T}(S_{\mathcal{T}}\cup\{S_i\})\subseteq \mathcal{T}$, then it contradicts the assumption that $S_{\mathcal{T}}$ is the maximum set containing all simple $A$-modules in $\mathcal{T}$. If $\mathcal{T}\subsetneq \mathscr{T}(S_{\mathcal{T}}\cup\{S_i\})$, then it contradicts the fact that $\mathscr{T}(S_{\mathcal{T}}\cup\{S_i\})\succ \mathscr{T}(S_{\mathcal{T}})$.

Hence we have
\[\mathscr{T}(S_{\mathcal{T}}) \leq \mathcal{T}\bigwedge \mathscr{T}(S_{\mathcal{T}}\cup\{S_i\})\lneq \mathscr{T}(S_{\mathcal{T}}\cup\{S_i\}),
\]
for any simple $A$-module $S_i \notin S_{\mathcal{T}}$. By the covering condition, we have
\[\mathcal{T}\bigwedge \mathscr{T}(S_{\mathcal{T}}\cup\{S_i\})= \mathscr{T}(S_{\mathcal{T}}).
\]
Then it follows from the completely semidistributivity of $\operatorname{tors} A$ that
\[\mathcal{T}\bigwedge (\bigvee_{S_{i} \notin S_{\mathcal{T}}}  \mathscr{T}(S_{\mathcal{T}}\cup\{S_i\}))= \mathscr{T}(S_{\mathcal{T}}).\]
But note that
\[\mathcal{T}\bigwedge (\bigvee_{S_{i} \notin S_{\mathcal{T}}}  \mathscr{T}(S_{\mathcal{T}}\cup\{S_i\}))= \mathcal{T}\bigwedge \operatorname{mod}A = \mathcal{T},
\]
then we obtain $\mathcal{T} = \mathscr{T}(S_{\mathcal{T}})$, which contradicts the assumption that $\mathscr{T}(S_{\mathcal{T}}) \subsetneq \mathcal{T}$.

Therefore, we proved that $\mathcal{T}=\mathscr{T}(S_{\mathcal{T}})$ for any $\mathcal{T} \in \operatorname{tors}A$, where $S_{\mathcal{T}}$ is the set of simple $A$-modules it contains. In particular, for every $B\in \operatorname{brick} A$. we have $\mathscr{T}(B) =\bigvee \mathscr{T}(S_{i})$, where each $S_i$ is a simple $A$-module in $\mathscr{T}(B)$. On the other hand, note that $\mathscr{T}(B)$ is completely join irreducible from Lemma \ref{bricks=join irre}, then there is only one simple module $S_i$ such that $\mathscr{T}(B) = \mathscr{T}(S_i)$. Hence we have $B\simeq S_i$.

$(b)\Rightarrow (c)$: It follows from Lemma \ref{every torsion and its bricks} that $\operatorname{tors}A$ is the Boolean lattice of $\{S_1, \dots, S_n\}$ partially ordered by inclusion of subsets. Then the result follows.
\end{prf} \medskip

Next we give the description of the algebras with such properties. Before this, we recall some notions.

Let $U$ be an $A$-module. Denote by $\operatorname{add} U$ the smallest full subcategory of $\operatorname{mod} A$ containing $U$, that is, the full subcategory of $\operatorname{mod} A$ whose objects are the finite direct sums of direct summands of $U$. Let $U'\in \operatorname{add} U$ and let $f: M \to U'$ be a morphism. The morphism $f$ is called a {\it left $\operatorname{add} U$-approximation} of $M$ if the induced morphism $\operatorname{Hom}_A(f, U^{''})$
is surjective for any $ U^{''}\in \operatorname{add} U$. The left $\operatorname{add} U$-approximation is called {\it minimal} if every $h \in \operatorname{End} U'$ such that $hf=f$ is an automorphism.

\begin{lem} \cite[Theorem 4.1]{Demonet2} \label{injection between tau rigid and bricks}
Let $A$ be a finite-dimensional $K$-algebra. Then there is an injection $i\tau$-$rigid$ $A$ $\rightarrow \operatorname{brick} A$ sending $M$ to $M/\operatorname{rad_{End(M)}}M$.
Moreover, if $A$ is $\tau$-tilting finite, the map is a bijection.
\end{lem}

\begin{lem}\label{left mutation of tautilting}\cite[Theorem 2.30 (b)]{Adachi}
Let $T= X\oplus U$ be a basic $\tau$-tilting $A$-module that is the Bongartz completion of $U$, where $X$ is indecomposable. Let $X\stackrel{f}\rightarrow U^{'} \rightarrow Y \rightarrow 0$ be an exact sequence, where $f$ is a minimal left $(\operatorname{add} U)$-approximation. If $U$ is sincere, then $Y$ is a direct sum of copies of an indecomposable $A$-module $Y_1$ and is not in $\operatorname{add} T$. In this case, $Y_1\oplus U$ is a basic $\tau$-tilting $A$-module.
\end{lem}

\begin{prop}\label{bricks simple local}
Let $A$ be a basic finite-dimensional $K$-algebra. Then the following statements are equivalent.
\begin{description}
\item $(a)$~Every brick module in $\operatorname{mod} A$ is simple.
\item $(b)$~The regular module $A$ is the unique basic $\tau$- tilting module.
\item $(c)$~The algebra $A$ is the finite direct product of finite-dimensional local algebras.
\end{description}
\end{prop}

\begin{prf}
$(a)\Rightarrow (b)$: Since all bricks in mod$A$ are simple, we know from Lemma \ref{every torsion and its bricks} that every torsion class is determined by simple modules, leading to $A$ is $\tau$-tilting finite. By Lemma \ref{injection between tau rigid and bricks}, the set of non-isomorphic $\tau$-rigid $A$-modules are precisely the set of non-isomorphic projective $A$-modules. Hence $A_A$ as the regular $A$-module is the unique basic $\tau$-tilting module.

$(b)\Rightarrow (c)$: Note that $A$ is basic. Let $A_A= P_1 \oplus \dots \oplus P_n$ be an indecomposable decomposition of $A$. Suppose that $A$ is not the product of some finite-dimensional local algebras, then $n > 1$. Without loss of generality, we may assume that $P_1$ satisfies the condition $\operatorname{Hom}_{A}(P_1, P_2 \oplus \dots \oplus P_n)\neq 0$. We write $Q= P_2 \oplus \dots \oplus P_n$ simply. Let $f: P_1\rightarrow Q^{'}$ be a minimal left $(\operatorname{add} Q)$-approximation of $P_1$. It is clear that $f\neq 0$ and $Q$ is sincere. Then we have an exact sequence of $\operatorname{mod} A$:
\[P_1 \stackrel{f}\rightarrow Q^{'}\stackrel{g}\rightarrow \operatorname{Coker} f \rightarrow 0.
\]
It follows from Lemma \ref{left mutation of tautilting} that there exists an indecomposable $Y_1$ in $\operatorname{mod} A$ such that $\operatorname{Coker} f= Y_{1}^{m}$ for some $m$, and $Q\oplus Y_1$ is a basic $\tau$-tilting $A$-module. Note that $A$ is the unique $\tau$-tilting $A$-module, we get that $Y_1\cong P_1$ and $\operatorname{Coker} f \cong P_1^{m}$. It follows that  $g$ is a retraction, thus, $P_{1}^{m}$ is a non-zero direct summand of $Q^{'}$. So we have $P_1 \in \operatorname{add}~Q$, which is a contradiction.

$(c)\Rightarrow (a)$: Note that $A$ is the finite direct product of finite-dimensional local algebras. For every indecomposable module $M \in \operatorname{mod}A$, $M$ is filtered by only one simple $A$-module. Particularly for any brick $B$, $B$ is filtered by its top $S$. Then we have $\mathscr{T}(B)= \mathscr{T}(S)$, that is, $B \simeq S$ by Lemma \ref{bricks=join irre}.
\end{prf}\medskip

Now we prove Theorem \ref{Main theorem}.
\begin{proof}[Proof of Theorem \ref{Main theorem}]
It follows from Proposition \ref{equivalent conditions of torsA} that $\operatorname{tors}A$ is distributive if and only if all bricks in $\operatorname{mod} A$ are simple modules. When $A$ is basic, then we obtain from Proposition \ref{bricks simple local} that they are equivalent to the condition that $A$ is the finite direct product of finite-dimensional local algebras.
\end{proof} \medskip

The lattice of torsion classes is a classical example of completely semidistributive lattices. The following corollary tells us that there are lots of semidistributive lattices which are not the torsion lattices.

\begin{cor} Let $L$ be a lattice which is distributive but not Boolean, then there is no finite-dimensional algebra $A$ such that $\operatorname{tors}A \simeq L$.
\end{cor}

Next, we investigate the distributivity of the poset $\operatorname{f-tors} A$ when it forms a lattice.  Let us recall some results of functorially finite torsion classes and $\tau$-tilting theory.

\begin{lem}\label{arrow separated of stau tilting}\cite[Theorem 2.35]{Adachi}
Let $U$ and $T$ be basic support $\tau$-tilting $A$-modules such that $U > T$ (i.e., $\operatorname{Fac}(U)\supsetneq \operatorname{Fac}(T)$). Then:
\begin{description}
\item (1)~there exists a right mutation $V$ of $T$ such that $\operatorname{Fac}(U)\supseteq \operatorname{Fac}(V)\supsetneq \operatorname{Fac}(T)$;
\item (2)~there exists a left mutation $V^{'}$ of $U$ such that $\operatorname{Fac}(U)\supsetneq \operatorname{Fac}(V^{'})\supseteq \operatorname{Fac}(T)$.
\end{description}
\end{lem}

The following lemma shows that $\operatorname{f-tors} A$ forms a lattice (respectively, complete lattice) ordered by inclusion if and only if it is a sublattice (respectively, complete sublattice) of $\operatorname{tors} A$.

\begin{lem}\label{Meet Join of ftorsA}\cite[Theorem 2.12]{Iyama}
Let $A$ be a finite-dimensional $K$-algebra.
\begin{description}
\item (1)~A subset $\{\mathcal{T}_i\mid i\in I\}$ of $\operatorname{f}$-$\operatorname{tors} A$ has a meet if and only if $\bigcap_{i\in I}\mathcal{T}_{i}$ is functorially finite.
\item (2)~A subset $\{\mathcal{T}_i\mid i\in I\}$ of $\operatorname{f}$-$\operatorname{tors} A$ has a join if and only if $^{\perp}(\bigcap_{i\in I}\mathcal{T}_{i}^{\perp})$ is functorially finite.
\end{description}
\end{lem}

For the poset $\operatorname{f-tors} A$, we have the following.
\begin{lem} \label{ftorsA semimodular then finie}
Let $A$ be a finite-dimensional $K$-algebra. If the poset $\operatorname{f-tors} A$ forms a lower semimodular lattice, then $A$ is $\tau$-tilting finite.
\end{lem}
\begin{prf}
Since $\operatorname{f-tors} A$ is a lattice, it follows from Lemma \ref{Meet Join of ftorsA} that $\operatorname{f-tors} A$ is a sublattice of $\operatorname{tors} A$. Thus $\operatorname{f-tors} A$ is also completely semidistributive.
Note that $A$ as $A$-module is a $\tau$-tilting module. The corresponding functorially finite torsion class $\operatorname{mod} A$ is the maximum of $\operatorname{f-tors}A$. By Lemma \ref{Hasse regular of ftorsA}, there are only $n$ arrows starting from $\operatorname{mod} A$, where $n$ is the number of non-isomorphic simple $A$-modules. Denote by $X_1, \dots, X_n$ all the functorially finite torsion classes adjacent to $\operatorname{mod} A$, respectively. It is clear that
$X_i \wedge X_j \neq X_i \wedge X_k$, for any $i, j, k \in \{1, \dots, n\}$ and $j \neq k$. Otherwise, it follows from $X_i \wedge X_j = X_i \wedge X_k$ and the complete semidistributivity that
\[X_i \wedge X_j = X_i \wedge (X_j \vee X_k)= X_i \wedge \operatorname{mod}A = X_i,
\]
which contradicts the fact that $X_i$ and $X_j$ are not comparable. Moreover, we also get that $X_i \wedge X_j$ and $X_i \wedge X_k$ are not comparable. Otherwise, we assume that $X_i \wedge X_j \lneq X_i \wedge X_k$.
Then it follows
\[ X_i \wedge X_j \lneq X_i \wedge X_k \lneq X_i,
\]
thus, $X_i \wedge X_j \nprec X_i$. On the other hand, since
\[X_i, X_j \prec X_i \vee X_j = \operatorname{mod}A,
\]
the lower semimodularity of $\operatorname{f-tors} A$ implies that $X_i \wedge X_j \prec X_i$, a contradiction. Hence $X_i \wedge X_j$ and $X_i \wedge X_k$ are not comparable.
Note that $\operatorname{f-tors}A$ is Hasse-regular. Then all adjacent elements of $X_i$ in the Hasse quiver are $\operatorname{mod} A$ and $X_i \wedge X_j$ for any $i \neq j$.

Next, we consider the set of adjacent elements of $X_i \wedge X_j$ for any $i \neq j$. Again by the lower semimodularity and complete semidistributivity, the adjacent elements of $X_i \wedge X_j$ are $X_i$, $X_j$, $X_i \wedge X_j \wedge X_k$ for any $k \neq i, j$. Continue these steps, we finally get that all adjacent elements of $\bigwedge_{i\in \{1, \dots, n\}} \{X_i\}$ are $\overline{X_i}$, for all $i\in \{1, \dots, n\}$, where $\overline{X_i}$ is the meet of $\{X_1, \dots, X_{i-1}, X_{i+1}, \dots, X_n\}$. Note that $\bigwedge_{i\in \{1, \dots, n\}} \{X_i\}$ is smaller than its all adjacent elements, and $\operatorname{f-tors} A$ is always Hasse-regular. It follows from Lemma \ref{arrow separated of stau tilting} that there is no elements strict smaller than $\bigwedge_{i\in \{1, \dots, n\}} \{X_i\}$, that is, $\bigwedge_{i\in \{1, \dots, n\}} \{X_i\}= 0$.

In fact, all elements in $\operatorname{f-tors}A$ are listed completely in this process. Indeed, for any $\mathcal{T} \in \operatorname{f-tors}A$ with $\mathcal{T}\neq \operatorname{mod} A$, then we know from Lemma \ref{arrow separated of stau tilting} again that $\mathcal{T} \leq X_i$ for some $i\in \{1,\dots, n\}$. If $\mathcal{T} \lneq X_i$, then we get that $\mathcal{T} \leq X_i \wedge X_j$ for another $j \in \{1, ..., n\}$. By induction, we proved that $\mathcal{T} = X_i \wedge X_j \wedge \dots \wedge X_k$ for some $i, j, \dots, k \in \{1, , n\}$. Certainly, $\operatorname{f-tors}A$ is finite, which implies that $A$ is $\tau$-tilting finite.
\end{prf}\medskip

Immediately, we get the following corollary.
\begin{cor} Let $A$ be a basic finite-dimensional $K$-algebra. Then $\operatorname{f-tors} A$ is a distributive lattice if and only if $A$ is a finite direct product of finite-dimensional local algebras.
\end{cor}
\begin{prf} Note that $A$ is $\tau$-tilting finite if and only if $\operatorname{tors} A = \operatorname{f-tors} A$. Then the result follows from Theorem \ref{Main theorem} and Lemma \ref{ftorsA semimodular then finie}
\end{prf} \medskip

\noindent {\bf Acknowledgement} The first author would like to thank Baptiste Rognerud for his valuable suggestions. Both authors also express deep gratitude to the reviewers for careful reading and insightful comments.

\noindent {\bf Declaration of interests} The authors have no competing interests to declare that are relevant to the content of this article.

\end{document}